\documentclass{amsart}
\usepackage{amsmath}
\usepackage{amsthm}
\usepackage{amsfonts}
\usepackage{amssymb}
\usepackage{amsbsy}
\usepackage{amscd}
\usepackage{eucal}
\usepackage[dvips]{graphicx}
\usepackage{amsaddr}

\newtheorem{theorem}{Theorem}
\newtheorem*{theorem2prime}{Theorem 2'}

\newtheorem{proposition}[theorem]{Proposition}

\theoremstyle{definition}

\numberwithin{equation}{section}

\newcommand{\beq}{\begin{equation}}
\newcommand{\eeq}{\end{equation}}

\newcommand{\rmd}{\mathrm{d}}
\newcommand{\rmi}{\mathrm{i}}

\newcommand{\N}{\mathbb{N}}
\newcommand{\Z}{\mathbb{Z}}
\newcommand{\R}{\mathbb{R}}

\newcommand{\hf}{\widehat{f}}

\newcommand{\ud}{\frac{1}{2}}

\newcommand{\sku}{\vspace*{0.1cm}}
\newcommand{\skd}{\vspace*{0.2cm}}
\newcommand{\skt}{\vspace*{0.3cm}}

\begin{document}

\title[Ultraspherical coefficients]{The expansion in ultraspherical polynomials:
a simple procedure for the fast computation of the ultraspherical coefficients}

\author[E. De Micheli]{Enrico De Micheli}
\address{\sl Consiglio Nazionale delle Ricerche \\ Via De Marini, 6 - 16149 Genova, Italy \\
E-mail: enrico.demicheli@cnr.it}
\author[G. A. Viano]{Giovanni Alberto Viano}
\address{\sl Dipartimento di Fisica -- Universit\`a di Genova,\\
Istituto Nazionale di Fisica Nucleare -- Sezione di Genova, \\
Via Dodecaneso, 33 - 16146 Genova, Italy \\
E--mail: viano@ge.infn.it}


\begin{abstract}
We present a simple and fast algorithm for the computation of the coefficients of
the expansion of a function $f(\cos u)$ in ultraspherical (Gegenbauer) polynomials. We prove that these
coefficients coincide with the Fourier coefficients of an Abel--type transform
of the function $f(\cos u)$. This allows us to fully exploit the computational efficiency of 
the Fast Fourier Transform, computing the first $N$ ultraspherical coefficients
in just $\mathcal{O} (N\log_2 N)$ operations.
\end{abstract}

\maketitle

\section{Introduction}
\label{se:introduction}
In this paper we investigate the generalization to ultraspherical polynomials 
(also known as Gegenbauer polynomials) of the results obtained in a previous paper \cite{DeMicheli} 
regarding the efficient computation of the coefficients of Legendre expansions.  

Ultraspherical expansions play a relevant role in various subjects of applied and computational
mathematics. These expansions have been used successfully for the solution of linear \cite{Elliott} 
and nonlinear \cite{Denman1,Denman2,Sinha} differential 
equations, of integral equation \cite{Srivastava}, and in spectral methods for 
partial differential equations \cite{Ben,Vozovoi}. Gegenbauer filtering (for the suppression of the Gibbs phenomenon)
\cite{Gelb,Gottlieb,Gottlieb2} has been proved to provide an exponentially convergent approximation 
(in the \emph{maximum} norm) of a piecewise analytic function, starting from the Fourier partial sum of the function itself.
The related Gegenbauer reconstruction method has found natural application in the image segmentation problem, in particular
of MRI images \cite{Archibald,Huang}.

One of the limits of the application of ultraspherical expansions is the high cost of computing the expansion
coefficients, this question becoming particularly critical in the case of multivariate functions.

In this paper we give an efficient procedure for computing these coefficients. We first obtain a Dirichlet--Murphy--type
integral representation of the ultraspherical polynomials $P_n^{\,(d)}(\cos u)$ of degree $n$ and order $d$.
Then we prove that the coefficients of the ultraspherical expansion of a function $f(\cos u)$ ($u\in[0,\pi]$) coincide with
the Fourier coefficients (restricted to nonnegative index) of an Abel--type transform of the function $f$.

These results produce straightforwardly an algorithm for the computation of ultraspherical coefficients 
which is very simple and very fast. 
In fact, the first $N$ ultraspherical coefficients are obtained in only $\mathcal{O}(N\log_2 N)$ operations 
by a single Fast Fourier Transform of the Abel--type integral function, 
the latter being easily computable by standard quadrature techniques \cite{Monegato}. 

Finally, in the Appendix the results obtained for the ultraspherical polynomials 
are given also in terms of Gegenbauer polynomials $C_n^{\,(\lambda)}(x)$ with $x\in[-1,1]$.

\section{Connection between Ultraspherical expansions and Fourier series}
\label{se:dimension}

In dimension $d$, the expansion in ultraspherical polynomials of a function $f=f(\cos u)$ ($u\in[0,\pi]$) reads:
\beq
f(\cos u) = \frac{1}{2^{(d-2)}\,\pi^{\frac{(d-1)}{2}}\Gamma(\frac{d-1}{2})}
\sum_{n=0}^\infty\frac{\Gamma(n+d-2)}{n!}\left(n+\frac{d-2}{2}\right)
\,a_n^{(d)}\,P_n^{\,(d)}(\cos u),
\label{d.-1}
\eeq
the coefficients $a_n^{(d)}$ being defined by \cite{Bros}: 
\beq
a_n^{(d)} =\frac{2\pi^{\frac{(d-1)}{2}}}{\Gamma(\frac{d-1}{2})}
\int_0^\pi f(\cos u) \ P_n^{\,(d)}(\cos u)\,(\sin u)^{(d-2)}\,\rmd u
\qquad (n\in\N_0),
\label{d.0}
\eeq
where $P_n^{\,(d)}(\cos u)$ denotes the ultraspherical polynomial of degree $n$ and order $d$, 
which is defined by the integral representation \cite{Faraut,Vilenkin}:
\beq
P_n^{\,(d)}(\cos u)=\frac{\Gamma(\frac{d-1}{2})}{\sqrt{\pi}\,\Gamma(\frac{d-2}{2})}
\int_0^\pi (\cos u+\rmi\sin u\cos\eta)^n\,(\sin\eta)^{(d-3)}\,\rmd\eta
\quad (n\in\N_0).
\label{d.1}
\eeq
Our goal now is to show that the coefficients $a_n^{(d)}$ are the Fourier coefficients of a suitable
Abel--type transform of $f$. To this end, we first prove the following proposition (see also \cite{Bros}).

\begin{proposition}
\label{pro:1}

The following integral representation of the ultraspherical polynomials $P_n^{\,(d)}(\cos u)$ holds:
\beq
P_n^{\,(d)}(\cos u) = \frac{(-\rmi)^{(d-2)}}{\sqrt{\pi}}\frac{\Gamma(\frac{d-1}{2})}{\Gamma(\frac{d-2}{2})}
\frac{1}{(\sin u)^{(d-3)}}
\int_u^{2\pi-u} \!\! e^{\rmi(n+\frac{d-2}{2})t} \ [2(\cos u-\cos t)]^{\frac{d-4}{2}}\,\rmd t.
\label{d.22}
\eeq 
\end{proposition}

\begin{figure}[tb]
\begin{center}
\leavevmode
\includegraphics[width=6cm]{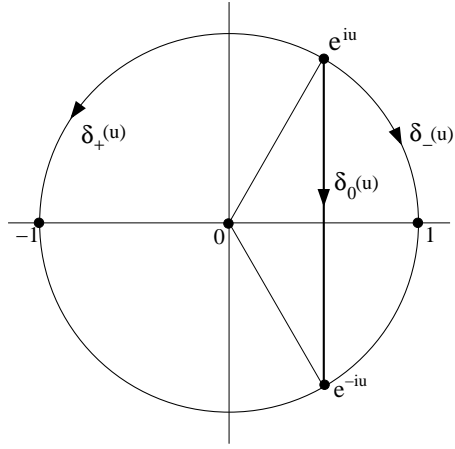}
\caption{\label{fig:1} \small Integration path for evaluating the integral representation (\protect\ref{d.22}) of
the ultraspherical polynomials.}
\end{center}
\end{figure}

\begin{proof}
In the integral representation \eqref{d.1} of $P_n^{\,(d)}(\cos u)$ substitute 
to $\eta$ the complex integration variable $\tau$ defined by
\beq
e^{\rmi\tau}=\cos u + \rmi\sin u\cos\eta.
\label{d.p1}
\eeq
It can be checked that
\beq
2e^{\rmi\tau}(\cos\tau-\cos u)=(e^{\rmi\tau}-e^{\rmi u})(e^{\rmi\tau}-e^{-\rmi u})=\sin^2\! u \, \sin^2\!\eta.
\label{d.p2}
\eeq
Now, since $e^{\rmi\tau}\rmd\tau=-\sin u\sin\eta\,\rmd\eta$, the integrand on the r.h.s. of Eq. \eqref{d.1} 
can be written as follows:
\beq
\begin{split}
& -(\sin u)^{-(d-3)}e^{\rmi(n+1)\tau}\left[(e^{\rmi\tau}-e^{\rmi u})
(e^{\rmi\tau}-e^{-\rmi u})\right]^\frac{d-4}{2}\,\rmd\tau \\
&\quad =  -(\sin u)^{-(d-3)}e^{\rmi(n+\frac{d-2}{2})\tau}\left[2(\cos\tau-\cos u)\right]^\frac{d-4}{2}\,\rmd\tau.
\end{split}
\label{d.p3}
\eeq
In order to determine the integration path, consider an intermediate step where $e^{\rmi\tau}$ 
is chosen as the integration variable; the original path (corresponding to $\eta\in[0,\pi]$) is the 
(oriented) linear segment $\delta_0(u)$ starting at $e^{\rmi u}$ and ending at $e^{-\rmi u}$. 
Since (as shown by \eqref{d.p3}) the integrand is an analytic function of $e^{\rmi\tau}$ in the 
disk $|e^{\rmi\tau}|<1$ (since $n\in\N$), the integration path $\delta_0(u)$ can be replaced by the circular path 
$\delta_+(u)=\{e^{\rmi\tau};\tau=t,u\leqslant t\leqslant 2\pi-u\}$ (see Fig. \ref{fig:1}). Moreover, by using the fact that
$[(e^{\rmi\tau}-e^{\rmi u})(e^{\rmi\tau}-e^{-\rmi u})]^\frac{d-4}{2}$ is positive for $e^{\rmi\tau}\in\delta_0(u)\cup\R$
and therefore at $e^{\rmi\tau}=e^{\rmi\pi}$, we conclude from the left equality in \eqref{d.p2} that in the r.h.s.
of \eqref{d.p3} the following specification holds (for $\tau=t; u\leqslant t\leqslant 2\pi-u$):
\beq
[2(\cos t-\cos u)]^\frac{d-4}{2}=(-\rmi)^{d-4}[2(\cos u-\cos t)]^\frac{d-4}{2}.
\label{d.p4}
\eeq
Finally, by taking into account the latter expression, the integral representation \eqref{d.1} can then be replaced by the 
integral representation \eqref{d.22}.
\end{proof}

We can now prove the following theorem (see also \cite{Bros}).

\begin{theorem}
\label{the:1}

The ultraspherical coefficients $\{a_n^{(d)}\}_{n=0}^\infty$ (see \eqref{d.0}) coincide with 
the Fourier coefficients (restricted to the nonnegative integer) of the form:
\beq
a_n^{(d)} = \int_{-\pi}^\pi \hf^{\,(d)}(t) \, e^{\rmi nt}\,\rmd t \qquad (n \in \N_0),
\label{d.3}
\eeq
where:
\beq
\hf^{\,(d)}(t)=\left[-\rmi\,\varepsilon(t)\right]^{d-2}\frac{2\pi^{\frac{d-2}{2}}}{\Gamma(\frac{d-2}{2})}
\ e^{\,\rmi\left(\frac{d-2}{2}\right)t}
\! \int_0^t f(\cos u)\,\left[2\left(\cos u-\cos t\right)\right]^\frac{d-4}{2}\sin u\,\rmd u,
\label{d.4}
\eeq
$\varepsilon(t)$ being the sign function.
\end{theorem}

\begin{proof}
By plugging representation \eqref{d.22} into formula \eqref{d.0} we have:
\beq
\label{d.5}
\begin{split}
a_n^{(d)} & =(-\rmi)^{d-2}\,\frac{2\pi^{\frac{d-2}{2}}}{\Gamma(\frac{d-2}{2})} \int_0^\pi\!\! \rmd u \, f(\cos u)\sin u
\int_u^{2\pi-u}\!\!e^{\rmi(n+\frac{d-2}{2})t}\ [2(\cos u-\cos t)]^{\frac{d-4}{2}}\,\rmd t \\
&=(-\rmi)^{d-2}\,\frac{2\pi^{\frac{d-2}{2}}}{\Gamma(\frac{d-2}{2})}\left\{
\int_0^\pi\!\rmd t \ e^{\rmi(n+\frac{d-2}{2})t}\!\!\int_0^t\!\rmd u \ f(\cos u)\sin u \
[2(\cos u-\cos t)]^{\frac{d-4}{2}} \right. \\
&\left.\quad +\int_\pi^{2\pi}\!\rmd t \ e^{\rmi(n+\frac{d-2}{2})t}\!\! 
\int_0^{2\pi-t}\!\!\rmd u \ f(\cos u)\sin u \ [2(\cos u-\cos t)]^{\frac{d-4}{2}}\right\}.
\end{split}
\eeq
Now, changing the variables $(t,u)\to (t+2\pi,u)$ in the second term inside the parentheses 
on the rightmost side of \eqref{d.5}, the latter becomes: 
\beq
e^{\rmi\pi(d-2)}\!\int_{-\pi}^0\rmd t \ e^{\rmi(n+\frac{d-2}{2})t}\!\int_0^{-t}\!\rmd u \ f(\cos u)\sin u \
[2(\cos u-\cos t)]^{\frac{d-4}{2}}.
\label{d.6}
\eeq
Finally, from \eqref{d.5} and \eqref{d.6} we have:
\beq
\label{d.7}
\begin{split}
a_n^{(d)} & =(-\rmi)^{d-2}\,\frac{2\pi^{\frac{d-2}{2}}}{\Gamma(\frac{d-2}{2})}\left\{
\int_0^\pi\!\rmd t \ e^{\rmi n t}\left[e^{\rmi\frac{d-2}{2}t}\!\int_0^t\rmd u \
f(\cos u)\sin u \ [2(\cos u-\cos t)]^{\frac{d-4}{2}}\right] \right. \\
&\left.\qquad +e^{\rmi\pi(d-2)}\!\int_{-\pi}^0\!\rmd t \ e^{\rmi nt}
\left[e^{\rmi\frac{d-2}{2}t}\!\int_0^{t}\rmd u \ f(\cos u)\sin u \ [2(\cos u-\cos t)]^{\frac{d-4}{2}}\right]\right\} \\
& =\int_{-\pi}^\pi \hf^{\,(d)}(t)\,e^{\rmi n t}\,\rmd t,
\end{split}
\eeq
with $\hf^{\,(d)}(t)$ given by \eqref{d.4}.
\end{proof}

\noindent
It is easy to check that the $2\pi$--periodic $\hf^{\,(d)}(t)$ function enjoys the following symmetry properties:
\beq
\hf^{\,(d)}(t)= (-1)^d\,e^{\,\rmi(d-2)t}\, \hf^{\,(d)}(-t) \qquad (t\in\R),
\label{2.12bis}
\eeq
and, consequently,
\beq
a_n^{(d)} = (-1)^d\,a_{-(n+d-2)}^{(d)} \qquad (n\in\Z).
\label{2.13bis}
\eeq

\newpage

\appendix

\section{From ultraspherical to Gegenbauer polynomials}
\label{app:a}

\subsection{The Gegenbauer polynomials $\boldsymbol{C_n^{(\lambda)}}$}
\label{subse:gegenbauer}

For the convenience of the reader we summarize hereafter the main properties of the Gegenbauer 
polynomials $C_n^{(\lambda)}(x)$, which are strictly related to the ultraspherical 
polynomials $P_n^{\,(d)}(x)$ that we introduced in Section \ref{se:dimension}. 
The following formulae are listed in Ref. \cite{Wolfram} and can also be 
found in refs. \cite[Chapter 22]{Abramowitz} and \cite{Elliott,Gottlieb,Szego}.

\skd

\noindent
The Gegenbauer polynomials of order $\lambda$ can be defined in terms of their generating function:
\beq
\sum_{n=0}^\infty C_n^{\,(\lambda)}(x)\,t^n = (1-2xt+t^2)^{-\lambda}.
\label{a.1}
\eeq
From \eqref{a.1} we see that the Legendre polynomials $P_n(x)$ are the particular case of Gegenbauer polynomials
with $\lambda=\ud$, i.e.: $C_n^{\,(\ud)}(x)=P_n(x)$. 
\skd

\noindent
They satisfy the recurrence relation:
\beq
\begin{split}
& n\,C_n^{\,(\lambda)}(x)= 2(n+\lambda-1)\,x\,C_{n-1}^{\,(\lambda)}(x) -(n+2\lambda-2)C_{n-2}^{\,(\lambda)}(x), \\
& C_0^{\,(\lambda)}(x)=1, \quad C_1^{\,(\lambda)}(x)=2\lambda\, x.
\end{split}
\label{a.2}
\eeq

\skd

\noindent
In terms of Gaussian hypergeometric function they can be written as:
\beq
C_n^{\,(\lambda)}(x)=\frac{2^{1-2\lambda}\sqrt{\pi}\,\Gamma(n+2\lambda)}{n!\,\Gamma(\lambda)} 
\ _2F_1\!\left(-n,2\lambda+n;\lambda+\ud;\frac{1-x}{2}\right).
\label{a.4}
\eeq

\skd

\noindent
They can be written explicitly as
\beq
C_n^{\,(\lambda)}(x)=\sum_{k=0}^{\lfloor n/2 \rfloor} (-1)^k 
\frac{\Gamma(\lambda+n-k)}{k!\,(n-2k)!\,\Gamma(\lambda)}\,(2x)^{n-2k},
\label{a.5}
\eeq
and have the following integral representation:
\beq
C_n^{\,(\lambda)}(\cos u) = \frac{2^{1-2\lambda}\Gamma(n+2\lambda)}{n!\,[\Gamma(\lambda)]^2}
\int_0^\pi\left(\cos u+\rmi\sin u\cos\eta\right)^n\,(\sin\eta)^{2\lambda-1}\,\rmd\eta.
\label{a.5bis}
\eeq

\skd

\noindent
The Gegenbauer polynomials can be computed by the Rodrigues formula:
\beq
C_n^{\,(\lambda)}(x)=\frac{(-2)^n}{n!}\frac{\Gamma(n+\lambda)\,\Gamma(n+2\lambda)}{\Gamma(\lambda)\,\Gamma(2n+2\lambda)}
(1-x^2)^{-\lambda+\ud}\frac{\rmd^n}{\rmd x^n}\left[(1-x^2)^{n+\lambda-\ud}\right],
\label{a.7}
\eeq
which follows by induction from
\beq
\frac{\rmd}{\rmd x}C_n^{\,(\lambda)}(x)=2\lambda\,C_{n-1}^{\,(\lambda+1)}(x).
\label{a.6bis}
\eeq

\skd

\noindent
For fixed $\lambda$, 
the Gegenbauer polynomials are orthogonal on the interval $[-1,1]$ with respect to the weight function:
\beq
w^{(\lambda)}(x)=(1-x^2)^{\lambda-\ud},
\label{a.8}
\eeq
that is, for $n\neq  m$:
\beq
\int_{-1}^1 C_n^{\,(\lambda)}(x) \, C_m^{\,(\lambda)}(x)\, (1-x^2)^{\lambda-\ud} \,\rmd x = 0,
\label{a.9}
\eeq
and are normalized by:
\beq
\int_{-1}^1 \left[C_n^{\,(\lambda)}(x)\right]^2 \, (1-x^2)^{\lambda-\ud} \,\rmd x
= \frac{\pi \, 2^{1-2\lambda}\,\Gamma(n+2\lambda)}{n! \, (n+\lambda) \left[\Gamma(\lambda)\right]^2}.
\label{a.10}
\eeq

\skt

\subsection{Relation between Gegenbauer polynomials $\boldsymbol{C_n^{(\lambda)}}$ and ultraspherical
polynomials $\boldsymbol{P_n^{(d)}}$}

Comparing the integral representations \eqref{d.22} and \eqref{a.5bis}, it is easily seen that the relation between
ultraspherical $P_n^{\,(d)}(\cos u)$ and Gegenbauer polynomials $C_n^{\,(\lambda)}(\cos u)$ is: 
\begin{subequations}
\label{a.11}
\begin{align}
C_n^{\,(\lambda)}(\cos u) &= \frac{\Gamma(n+2\lambda)}{n!\,\Gamma(2\lambda)}\ P_n^{\,(2\lambda+2)}(\cos u),
\label{a.11a} \\[+5pt]
P_n^{\,(d)}(\cos u) &= \frac{n!\,\Gamma(d-2)}{\Gamma(n+d-2)}\ C_n^{\,(\frac{d-2}{2})}(\cos u).
\label{a.11b}
\end{align}
\end{subequations}

Then, rephrasing Proposition \ref{pro:1} in terms of Gegenbauer polynomials, from \eqref{d.22} and \eqref{a.11}
it follows that the polynomials $C_n^{\,(\lambda)}(\cos u)$ ($u\in[0,\pi]$) have the following integral representation:
\beq
C_n^{\,(\lambda)}(\cos u) =
\frac{\Gamma(2\lambda+n)\Gamma(\lambda+\ud)}{\sqrt{\pi}\,n!\,\Gamma(2\lambda)\,
\Gamma(\lambda)}\frac{(-\rmi)^{2\lambda}}{(\sin u)^{2\lambda-1}}
\int_u^{2\pi-u}\!\! e^{\rmi(n+\lambda)t}\,[2(\cos u-\cos t)]^{\lambda-1}\,\rmd t.
\label{a.14}
\eeq 

\sku 

\subsection{Expansions in Gegenbauer polynomials $\boldsymbol{C_n^{(\lambda)}}$}
\label{subse:expansion}

The Gegenbauer expansion of a function $f(x)$, defined in $x\in[-1,1]$, reads \cite{Gottlieb}:
\beq
f(x) = \sum_{n=0}^\infty f_n^{\,(\lambda)} C_n^{\,(\lambda)}(x)
\qquad (x\in[-1,1]; \lambda>0),
\label{a.15}
\eeq
where the Gegenbauer coefficients are given by 
\beq
f_n^{\,(\lambda)}=\frac{n!\,\Gamma(2\lambda)\,\Gamma(\lambda)\,(n+\lambda)}
{\sqrt{\pi}\,\Gamma(n+2\lambda)\,\Gamma(\lambda+\ud)}
\int_{-1}^{1} f(x)\,C_n^{\,(\lambda)}(x) \ (1-x^2)^{(\lambda-\ud)}\,\rmd x
\qquad (n\in\N_0).
\label{a.16}
\eeq
Now, from \eqref{a.16}, \eqref{a.11} and using definition \eqref{d.0} of the coefficients $a_n^{(d)}$
we see that the Gegenbauer coefficients are related to the ultraspherical coefficients by
\beq
f_n^{\,(\lambda)} = \frac{\Gamma(\lambda)\,(n+\lambda)}{2\,\pi^{(\lambda+1)}}\,a_n^{(2\lambda+2)}.
\label{a.18}
\eeq
We can now rewrite Theorem \ref{the:1} in terms of Gegenbauer coefficients $f_n^{\,(\lambda)}$.

\begin{theorem2prime}
\label{the:2'}

The coefficients $\{\phi_n^{(\lambda)}\}_{n=0}^\infty$, defined as 
$\phi_n^{(\lambda)}\doteq\frac{f_n^{\,(\lambda)}}{(n+\lambda)}$ (see \eqref{a.16}) coincide with the following
Fourier coefficients (restricted to the nonnegative integer): 
\beq
\phi_n^{(\lambda)}\doteq\frac{f_n^{\,(\lambda)}}{(n+\lambda)}
= \int_{-\pi}^\pi \widehat{\phi}^{\,(\lambda)}(t) \, e^{\rmi nt}\,\rmd t \qquad (n \in \N_0),
\label{a.19}
\eeq
where:
\beq
\widehat{\phi}^{\,(\lambda)}(t)
=\frac{[-\rmi\,\varepsilon(t)]^{2\lambda}}{\pi}\,e^{\rmi\lambda t}
\! \int_{\cos t}^1 f(x)\,\left[2\left(x-\cos t\right)\right]^{\lambda-1}\,\rmd x,
\label{a.20}
\eeq
$\varepsilon(t)$ being the sign function.
\end{theorem2prime}

\end{document}